\documentclass[10pt,twoside]{article}

\usepackage{amsbsy,amsfonts,amsmath,amssymb,eucal}
\usepackage[all]{xy}

\def\itemn#1{\item[\hspace{0.6mm} {\rm (#1)}]}
\def\qed{\hfill $\square$ \vspace{5mm}}

\def\str{{\rm str}}
\def\tir{\mbox{-}}
\def\zmod#1{\mathbb{Z}/#1\mathbb{Z}}

\def\longto{\longrightarrow}
\def\isomto{\stackrel{\sim}{\longto}}

\def\To{\Rightarrow}
\renewcommand{\tilde}{\widetilde}

\newtheorem{definition}{Definition}[subsection]

\newtheorem{remark}[definition]{Remark}

\newtheorem{remarks}[definition]{Remarks}

\newtheorem{example}[definition]{Example}

\newtheorem{examples}[definition]{Examples}

\newtheorem{nothing}[definition]{$\!\!$}

\newenvironment{proo}{{\flushleft \bf Proof :}}{\hfill $\square$ \vspace{5mm}}

\newtheorem{definition*}{Definition}[section]
\newenvironment{defi*}{\begin{definition*} \rm}{\end{definition*}}
\newtheorem{prop*}[definition*]{Proposition}
\newtheorem{lemm*}[definition*]{Lemma}
\newtheorem{coro*}[definition*]{Corollary}
\newtheorem{theo*}[definition*]{Theorem}
\newtheorem{remark*}[definition*]{Remark}
\newenvironment{rema*}{\begin{remark*} \rm}{\end{remark*}}
\newtheorem{remarks*}[definition*]{Remarks}
\newenvironment{remas*}{\begin{remarks*} \rm}{\end{remarks*}}
\newtheorem{example*}[definition*]{Example}
\newenvironment{exam*}{\begin{example*} \rm}{\end{example*}}
\newtheorem{examples*}[definition*]{Examples}
\newenvironment{exams*}{\begin{examples*} \rm}{\end{examples*}}

\DeclareMathOperator{\Hom}{Hom}
\DeclareMathOperator{\Isom}{Isom}
\DeclareMathOperator{\Aut}{Aut}
\DeclareMathOperator{\Id}{id}
\DeclareMathOperator{\Spec}{Spec}

\DeclareMathOperator{\Cat}{Cat}

\DeclareMathOperator{\GL}{GL}

\DeclareMathOperator{\pr}{pr}

 \def\clC{{\cal C}} \def\clD{{\cal D}} 
\def\clF{{\cal F}}  \def\clH{{\cal H}} 
\def\clM{{\cal M}} \def\clN{{\cal N}} \def\clO{{\cal O}} \def\clP{{\cal P}}

   \def\clX{{\cal X}}
 \def\clZ{{\cal Z}}

\def\fra{\mathfrak{a}} \def\frb{\mathfrak{b}} 
  
  \def\frS{\mathfrak{S}}
\def\frT{\mathfrak{T}} 

\def\frC{\mathfrak{C}} \def\frE{\mathfrak{E}} \def\frT{\mathfrak{T}}
\def\frGrpd{\mathfrak{Grpd}} \def\frSt{\mathfrak{St}}
\def\frAlgSt{\mathfrak{AlgSt}} \def\frCat{\mathfrak{Cat}}

\begin{document}

\begin{center}
{\bf \Large A note on group actions on algebraic stacks}
\par\vspace{10mm}\par
{\bf Matthieu Romagny}
\par\vspace{2mm}\par
\emph{Stockholms Universitet, Matematiska institutionen,
106 91 Stockholm, Sweden}

\emph{e-mail address: romagny@matematik.su.se}
\end{center}

\par\vspace{7mm}\par

{\def\thefootnote{\relax}
\footnote{ \hspace{-6.8mm}
}
}

\par\vspace{1mm}\par
\footnotesize
\noindent {\bf Abstract}: we give the basic definitions of group actions
on (algebraic) stacks, and prove the existence of fixed points and quotients
as (algebraic) stacks.
\normalsize

\par\vspace{2mm}\par

\section{Introduction}

For people using algebraic stacks, it is inevitable, one day or another,
to meet a group acting on it. Thus it is no surprise that group actions
on algebraic stacks must be studied. What is more surprising is that, up to
now, no systematic treatment has been done. On many occasions, such actions
appear in the literature, for example the action of a torus on a stack of
stable maps in \cite{Ko}, (3.2) and the action of the symmetric group $\frS_d$
on a stack of multisections in \cite{LMB}, (6.6). We can also mention the
increasing need of these concepts in (orbifold) Gromov-Witten theory.
But certainly the most natural instance of an action on an algebraic stack
comes when we consider a scheme $X$ with an action of an algebraic
group $G$. Indeed, then the normalizer of $G$ in $\Aut(X)$ (= the automorphism
group scheme of~$X$) acts on the quotient stack $[X/G]$~; it is in some sense
"what remains" of the symmetries of $X$.

Another natural example of group action on an algebraic stack comes from
Hurwitz stacks (this is the example that initially motivated this note).
These are stacks parameterizing covers of algebraic curves. If one
restricts attention to Galois covers of group $G$, then the automorphism group
of $G$ acts on the resulting stack by twisting the action of $G$.
If the curves have marked points then the symmetric group acts
on the stack by permutation of the points.
Recently these stacks have been
extensively studied. One should mention works of Bertin \cite{Be},
Ekedahl \cite{E}, Wewers \cite{We}, Abramovich-Corti-Vistoli \cite{ACV}.

The aim of this note is to give basic definitions of actions, and existence theorems
for fixed points and quotients. Here is an overview of both its contents
and organization. Section 2 is rather informal and prepares the definitions
related to actions on stacks in section 3. The "main course" is in section~4
where we prove algebraicity of fixed point stacks under a proper flat
group scheme, and of quotients under a separated flat group scheme (in this
section, groups are finitely presented). Having in mind the application to covers of
curves in mixed characteristic, where the framework of finite constant
groups is outdated by the last developments (such as works of Raynaud, Henrio,
Wewers, Sa\"\i di and others), we study particularly closely the
case of (arbitrary) finite flat group schemes.

\par\vspace{2mm}\par

\noindent \emph{Notations}.
In the article, a scheme $S$ is fixed once for all. Most schemes, spaces, and stacks
will be over $S$, and quite often the mention of $S$ (for instance in fibred products)
will be omitted. If $X$ is an $S$-scheme and $T\to S$ a base change, we often
write $X_T=X\times_S T$. We write categories in calligraphic letters (such as $\clC$)
and 2-categories in fraktur letters (such as $\frC$). A category fibred in groupoids
over $S$ is simply called a groupoid over $S$. In such a groupoid $\clM$, the functor
of isomorphisms between two objects $x,y\in\clM(T)$ is denoted $\Isom_T(x,y)$ or
$\Isom_{\clM_T}(x,y)$ if mention of $\clM$ is needed. An algebraic stack is
an algebraic stack in the sense of \cite{LMB}, def. 4.1.

\par\vspace{2mm}\par

\noindent \emph{Acknowledgments}. I thank B. Toen, L. Moret-Bailly and
A. Vistoli for various discussions and indications. I thank Stockholm
University for hospitality and the EAGER Network for financial support.

\section{Preliminaries}

We need first to recall some basics concerning diagrams in a 2-category. Loosely speaking,
a \emph{diagram} in a 2-category~$\frC$ is a set of objects, with a set of 1-morphisms
between certain pairs of objects, and a set of 2-morphisms between certain pairs of 1-morphisms
(understood, with same source and target). We will write $\clD=\{\clM,f,\alpha\}$
to indicate that $\clM$ ranges through the set of objects, $f$ ranges through the set of
1-morphisms and $\alpha$ ranges through the set of 2-morphisms of the diagram $\clD$.
Notice that, by \emph{the set of 1-morphisms of the diagram}, we mean a set which is
saturated under composition, i.e. including all possible compositions we can make
with the original 1-morphisms, and similarly with the 2-morphisms. We call \emph{circuit}
a pair of morphisms of $\clD$ with same source and target. A circuit \emph{commutes} if its
two morphisms coincide. Here is a first example~:
$$
\xymatrix@R=2.5pc{
\clM \ar[r]^f \ar[d]_h & \clM' \ar[r]^{f'} \ar[d]_{h'} & \clM'' \ar[d]^{h''} \\
\clN \ar[r]^g \ar@2[ru]^{\alpha} & \clN' \ar[r]^{g'} \ar@2[ru]^{\alpha'} & \clN''
}
$$
Here $(gh,h'f)$ is a circuit of 1-morphisms (or 1-cicuit), and $\alpha$ is a 2-morphism
between $gh$ and $h'f$, so that the 1-circuit commutes if $\alpha=\Id$. A second example
may be given by the same diagram, plus a 2-morphism $\alpha'':g'gh\To h''f'f$ attached
to the exterior rectangular 1-circuit~; then there is a 2-circuit
$$
\xymatrix@C=1pc{
& g'h'f \ar@2[rd]^{\alpha'f} & \\
g'gh \ar@2[ru]^{g'\alpha} \ar@2[rr]_{\alpha''} & & h''f'f
}
$$
If $\star$ denotes the composition of 2-morphisms, we will write this 2-circuit
$(\alpha'',\alpha'\star \alpha)$, though to be rigorous we should write
$(\alpha'',\alpha'f\star g'\alpha)$. A diagram in $\frC$ is said to be
\emph{2-commutative} if any of its 2-circuits commutes, i.e. we have
$\alpha''=\alpha'\star\alpha$ in the example above. Given a diagram in $\frC$,
by forgetting the 2-morphisms we get a diagram in the underlying category $\Cat(\frC)$
of $\frC$. Sometimes we refer to diagrams of $\frC$
as 2-diagrams, and diagrams of $\Cat(\frC)$ as 1-diagrams. So for example we will
say that a 2-diagram in $\frC$ is 1-commutative if the associated 1-diagram is
1-commutative.

Our main aim in the article is to discuss group actions on \emph{algebraic stacks}~;
of course this will just be an action on the underlying stack, and even on the underlying
groupoid. Thus we simplify the approach by looking first at the case of the 2-category
$\frGrpd/S$ of groupoids over $S$. In this particular 2-category, 2-commutativity
means "(1-)commutativity up to (given) isomorphisms".

Let $\clM$ be such a groupoid, and $G$ be a functor in groups over $S$. We denote by $m$
the multiplication of $G$, and by $e$, or sometimes simply~1, its unit section. As is natural,
one defines an action of $G$ on $\clM$ to be a morphism of groupoids $\mu:G\times \clM\to \clM$
satisfying the usual commutative diagrams concerning compatibility with respect to the unit
section of $G$ and to associativity of the multiplication. Using the more natural notion
of 2-commutativity in $\frGrpd/S$, we are led to 2-commutative diagrams
\begin{equation} \label{EQ1}
\xymatrix{
G\times G\times \clM \ar[r]^>>>>>{m\times \Id_{\clM}} \ar[d]_{\Id_G\!\times \mu} &
G\times \clM \ar[d]^{\mu} \\
G\times \clM \ar[r]^{\mu} \ar@2[ru]^{\alpha} & \clM
}
\qquad\qquad
\xymatrix{
G\times \clM \ar[r]^>>>>>{\mu} \ar@2+<7mm,-7mm>^>>>\fra & \clM \\
\clM \ar[u]^{e\times \Id_{\clM}} \ar[ru]_{\Id_{\clM}} &
}
\end{equation}
where $\alpha:\mu\circ (\Id\times \mu)\To \mu\circ (m\times \Id)$ and
$\fra:\mu\circ (e\times \Id)\To \Id$ are 2-isomorphisms. Let's make the usual
notational convention that, when we have an action of $G$ on $\clM$, and sections
$x\in\clM$ and $g\in G$ over a scheme $T$, we will write $g.x$ or $gx$, rather
than $\mu(g,x)$; for an arrow $\varphi:x\to y$ we will write $g.\varphi$ or
$g\varphi$ rather than $\mu(g,\varphi)$. So what the above 2-commutative
diagrams mean is that we are given isomorphisms in $\clM$, natural in $(g,h,x)$,
$$
\alpha_{g,h}^x:g.(h.x)\isomto (gh).x
\qquad\qquad\mbox{and}\qquad\qquad
\fra^x:1.x\isomto x
$$
Similarly, a morphism of  $G$-groupoids
$f:\clM\to \clN$ should be given by a 2-commutative diagram
\begin{equation} \label{EQ2}
\xymatrix{
G\times \clM \ar[r]^>>>>>{\mu} \ar[d]_{\Id_G\!\times f} &
\clM \ar[d]^{f} \\
G\times \clN \ar[r]^{\nu} \ar@2[ru]^{\sigma} & \clN
}
\end{equation}
where $\sigma$ can be written more concretely by
$$
\sigma_g^x:g.f(x)\isomto f(g.x)
$$

Furthermore, if we want a triple $(\mu,\alpha,\fra)$ as above to define an action,
the isomorphisms $\alpha_{g,h}^x$ should be compatible with associativity in $G$.
In the same vein, a pair $(f,\sigma)$ will give rise to a morphism if
$\sigma$ is compatible to $\alpha,\fra$ and the corresponding 2-isomorphisms
$\beta,\frb$ for $\clN$. We arrive at the following provisional
definition of an action.

\begin{defi*}
Let $\clM$ be a groupoid over $S$ and $G$ a functor in groups over $S$.
\begin{trivlist}
\itemn{i} 
An \emph{action} of $G$ on $\clM$ is a triple $(\mu,\alpha,\fra)$
where $\mu:G\times \clM\to \clM$ is a morphism of groupoids satisfying the
above 2-commutative diagrams (\ref{EQ1}), and such that for all $x$ and $g,h,k$ we have
$$
\alpha_{g,hk}^x\circ g.\alpha_{h,k}^x=\alpha_{gh,k}^x\circ \alpha_{g,h}^{k.x}
\qquad\qquad\mbox{and}\qquad\qquad
1.\fra^x=\alpha_{1,1}^x
$$
We also say that $\clM$ is a $G$-\emph{groupoid}. If $\alpha$
and $\fra$ are the identity 2-isomorphisms, we say that the action (or the
$G$-groupoid) is \emph{strict}.
Usually we simply note $\clM$ for $(\clM,\mu,\alpha,\fra)$.
\itemn{ii}
A \emph{morphism of $G$-groupoids} between
$(\clM,\mu,\alpha,\fra)$ and $(\clN,\nu,\beta,\frb)$ is a pair
$(f,\sigma)$ where $f:\clM\to \clN$ is a morphism of groupoids over $S$
satisfying the above 2-commutative diagram  (\ref{EQ2}), and such that for all $x$
and $g,h$ we have
$$
f(\alpha_{g,h}^x)\circ\sigma_g^{h.x}\circ g.\sigma_h^x=\sigma_{gh}^x\circ \beta_{g,h}^{f(x)}
\qquad\qquad\mbox{and}\qquad\qquad
f(\fra^x)\circ\sigma_1^x=\frb^{f(x)}
$$
\itemn{iii} An \emph{isomorphism of $G$-groupoids} is a morphism of
$G$-groupoids which is also an equivalence of categories fibred over $S$.
\end{trivlist}
\end{defi*}

\par\vspace{3mm}\par

The 2-category of $G$-groupoids over $S$ is denoted $G\tir\frGrpd/S$.
Before getting a headache trying to check if these are really the
good compatibilities, note that the key point in making these definitions
is the following remark. Once we make the definition of a $G$-groupoid,
we can recognize the data $(\clM,\mu,\alpha,\fra)$ as giving exactly what
is called a \emph{lax presheaf in groupoids $\clF$ over $\clC=B_0G$},
where $B_0G$ is the groupoid associated to $G$, i.e. the groupoid
whose fibre over a scheme $T/S$ has only one object, and morphisms the
elements of $G(T)$. This lax presheaf (see for instance \cite{Ho},
Appendix B) is described as follows~:
\begin{trivlist}
\itemn{i} To an object of $B_0G$, i.e. a scheme $T$ over $S$,
is associated the groupoid $\clF(T)=\clM(T)$.
\itemn{ii} To a morphism of $B_0G$, i.e. an element $g\in G(T)$,
is associated the functor $\mu(g^{-1},.):\clM(T)\to \clM(T)$
denoted $\mu_g$.
\itemn{iii} For each $g,h\in G(T)$ (= pair of composable arrows),
there is a natural transformation
$\mu_g\circ \mu_h\isomto \mu_{hg}$ given by $\alpha_{g^{-1},h^{-1}}$.
\end{trivlist}
So now, the definition of a morphism of  $G$-groupoids is just
a translation of the definition of a morphism of lax presheaves as we
find it in \cite{Ho}.
This link with lax presheaves also explains why in fact, given a $G$-groupoid
$\clM$, we will always be able to find an equivalent $G$-groupoid $\clM^\str$
such that the 2-isomorphisms $\alpha$ and $\fra$ are the identities.

\begin{prop*} \label{strictification}
There is a "strictification" functor $G\tir\frGrpd/S\to G\tir\frGrpd/S$
sending any $G$-groupoid to an isomorphic $G$-groupoid with strict action.
\end{prop*}

\begin{proo}
Let $\clM$ be a $G$-groupoid, and define a $G$-groupoid $\clM^\str$
in the following way~:
\begin{trivlist}
\itemn{i} the sections of $\clM^\str$ over a scheme $T$
are pairs $(g,x)$ with $g\in G(T)$, $x\in\clM(T)$,
\itemn{ii} the arrows in $\clM^\str$ bewteen $(g,x)$ and $(h,y)$ are arrows
$\varphi:x\to (g^{-1}h).y$ in $\clM(T)$,
\itemn{iii} composition of two arrows $\varphi:(g,x)\to (h,y)$
and $\psi:(h,y)\to (k,z)$ is given by
\end{trivlist}
$$
\xymatrix{
x \ar[r]^>>>>>{\varphi}
& (g^{-1}h).y \ar[rr]^(.4){(g^{-1}h).\psi}
& & (g^{-1}h).((h^{-1}k).z) \ar[rr]^>>>>>>>>>>>{\alpha_{g^{-1}h,h^{-1}k}^z}
& & (g^{-1}k).z
}
$$
There is a strict action of $G$ on $\clM^\str$~: an element $\gamma\in G(T)$
sends an object $(g,x)$ to $(\gamma g,x)$, and sends an arrow
$\varphi:x\to (g^{-1}h).y$ to the same arrow as a morphism 
between $(\gamma g,x)$ and $(\gamma h,y)$. Furthermore it is clear
that $\clM^\str$ is functorial in $\clM$.

It only remains to check that $\clM$ and $\clM^\str$ are isomorphic. We define
a morphism of groupoids $u:\clM^\str\to \clM$ by mapping an object $(g,x)$ to
$g.x$, and an arrow $(g,x)\to (h,y)$ represented by $\varphi:x\to (g^{-1}h).y$
to the composition $\alpha_{g,g^{-1}h}^y\circ (g.\varphi)$. Clearly, $u$ is a
$G$-morphism. Furthermore it is essentially surjective because any object $x$
in $\clM$ is isomorphic via $\fra^x$ to $1.x$. Finally it is straightforward
to see that it is fully faithful, so $u$ is an isomorphism.
\end{proo}

Now assume that the scheme $S$, viewed as the category of $S$-schemes,
is endowed with a Grothendieck topology. In practice, for us this will
be the fppf or \'etale topology. Then it is clear that an action of $G$
on a groupoid $\clM$ extends uniquely to an action on the associated stack
$\tilde{\clM}$. We could make, in the context of stacks over $S$, statements
similar to all ones in this section, taking associated stacks at the right
moments. (For example the groupoid $B_0G$ would be replaced by the stack
of $G$-torsors $BG$.)

In the next section we develop the basics of the theory of actions on stacks,
starting from the idea that we can restrict to considering strict actions.
This is made legitimate by proposition~\ref{strictification}. Note that the
theory could certainly be developped with general weak actions, at the cost
of substantial technical complications. This seems unnecessary~: the practice
will show that all constructions we wish to make will yield strict actions,
and if that were not the case, strictifying at the right place would
bring back to this context.

\section{Group actions on stacks}

Below is the minimal number of definitions that we will need for $G$-stacks.
There are two ways to present these concepts, according to whether we define
morphisms before commutative diagrams, or after them. We take the first option
and explain the second in a remark.

\begin{defi*} \label{defiActionStack}
Let $\clM$ be a stack over $S$ and $G$ a sheaf in groups over $S$.
Let $m$  denote the multiplication of $G$, and $e$ its unit section.
\begin{trivlist}
\itemn{i}
An \emph{action} of $G$ on $\clM$ is a morphism of stacks
$\mu:G\times \clM\to \clM$ with 1-commutative diagrams
$$
\xymatrix{
G\times G\times \clM \ar[r]^>>>>>{m\times \Id} \ar[d]_{\Id\times \mu} &
G\times \clM \ar[d]^{\mu} \\
G\times \clM \ar[r]^{\mu} & \clM
}
\qquad\qquad
\xymatrix{
G\times \clM \ar[r]^>>>>>{\mu} & \clM \\
\clM \ar[u]^{e\times \Id} \ar[ru]_{\Id}}
$$
We say that $(\clM,\mu)$ is a \emph{$G$-stack}.
\itemn{ii}
A \emph{1-morphism of $G$-stacks}, or \emph{1-$G$-morphism}, between
$(\clM,\mu)$ and $(\clN,\nu)$ is a pair $(f,\sigma)$ where $f:\clM\to \clN$
is a morphism of stacks with a 2-commutative diagram
$$
\xymatrix{
G\times \clM \ar[r]^>>>>>{\mu} \ar[d]_{\Id\times f} &
\clM \ar[d]^{f} \\
G\times \clN \ar[r]^{\nu} \ar@2[ru]^{\sigma} & \clN
}
$$
such that for all sections $x\in \clM(T)$ and $g,h\in G(T)$ over a scheme $T$,
the isomorphisms $\sigma_g^x:g.f(x)\simeq f(g.x)$ satisfy the cocycle relation
$\sigma_g^{h.x}\circ g.\sigma_h^x=\sigma_{gh}^x$.
Composition of 1-morphisms of $G$-stacks is defined in the obvious way, namely,
$(f_2,\sigma_2)\circ (f_1,\sigma_1)=(f_3,\sigma_3)$ where $f_3=f_2\circ f_1$ and
$\sigma_{3,g}^x=f_2(\sigma_{1,g}^x)\circ \sigma_{2,g}^{f(x)}$.
\itemn{iii}
A \emph{2-morphism of $G$-stacks}, or \emph{2-$G$-morphism}, between
1-morphisms $(f_1,\sigma_1)$ and $(f_2,\sigma_2)$ is a 2-morphism of stacks
$\tau:f_1\To f_2$ compatible with the $\sigma_i$ i.e. such that for all sections
$x\in \clM(T)$ and $g\in G(T)$ over a scheme~$T$, we have
$\sigma_{2,g}^x\circ g.\tau^x=\tau^{g.x}\circ \sigma_{1,g}^x$. In this way we have
defined a 2-category of $G$-stacks over $S$, which will be denoted by $G\tir\frSt/S$
or simply $G\tir\frSt$ if the base $S$ is clear. In particular, given two $G$-stacks
$\clM,\clN$ there is the stack $\clH om_{G\tir\frSt}(\clM,\clN)$ of 1-morphisms and
2-morphisms between them.
\itemn{iv}
Two actions $\mu,\mu'$ on the stack $\clM$ are said to
be \emph{equivalent} if there exists $\sigma$ such that the pair
$(\Id,\sigma)$ is a 1-$G$-morphism between $(\clM,\mu)$ and $(\clM,\mu')$.
\itemn{v}
An \emph{isomorphism of $G$-stacks} is a 1-$G$-morphism which is also an
equivalence of groupoids over $S$.
\end{trivlist}
\end{defi*}

\begin{rema*} \label{remaCommDiag}
We can restate definitions (ii) and (iii) using only diagrams in the 2-category
of stacks. Indeed, let $\clD=\{\clM,f,\alpha\}$ be a diagram of stacks where
$\clM$, $f$, $\alpha$ range through objects, 1-morphisms, and 2-morphisms of
$\clD$ respectively. Consider
\begin{trivlist}
\itemn{a}
the diagram
$G\times \clD:=\{G\times \clM \, , \,
\Id_G\times f\, ,\, \Id_{\Id_G}\times \alpha\}$,
\itemn{b}
the diagram
$G\times G\times \clD:=\{G\times G\times \clM\, ,\, \Id_{G\times G}\times f\, ,\,
\Id_{\Id_{G\times G}}\times \alpha\}$.
\end{trivlist}
Assume that the objects are in fact $G$-stacks $(\clM,\mu)$ and for any objects
$(\clM,\mu)$, $(\clN,\nu)$ and any 1-morphism $f:\clM\to \clN$ we are given a
2-isomorphism $\sigma:\nu\circ (\Id_G\times f)\To f\circ\mu$. Then we can form a new
diagram $G\times G\times \clD\to G\times \clD\to \clD$. Precisely, at the stage
$G\times \clD\to \clD$ the 1-morphisms are the $\mu$'s, the 2-morphisms are the
$\sigma$'s. At the stage $G\times G\times \clD\to G\times \clD$, the 1-morphisms
are the $(\Id_G\times\mu)$'s, the 2-morphisms are the $(\Id_G\times \sigma)$'s. Then
we can \emph{define} a \emph{2-commutative diagram of $G$-stacks} to be given by data
$(\clD,\{\mu\},\{\sigma\})$ such that $G\times G\times \clD\to G\times \clD\to \clD$
is a 2-commutative diagram of stacks.
In particular, we can now \emph{define} the notions of \emph{1-morphisms of $G$-stacks}
and \emph{2-morphisms between 1-morphisms of $G$-stacks}, by the 2-commutativity of the
following two elementary diagrams~:
$$
\xymatrix@C=1.2pc{
(\clM,\mu) \ar[rr]^{(f,\sigma)} & & (\clN,\nu)
}
\qquad\mbox{and}\qquad
\xymatrix@C=1.2pc{
(\clM,\mu) \ar@/^.8pc/[rr]^{(f_1,\sigma_1)} \ar@/_.8pc/[rr]_{(f_2,\sigma_2)}
\ar@{}[rr]|{\;\;\;\;\;\;\tau} & \Downarrow & (\clN,\nu)
}
$$
which means checking 2-commutativity of the following "prisms"
\par\vspace{2mm}\par
$$
\begin{array}[t]{c}
\\
\xymatrix@C=1pc{
G\times G\times \clM \ar[rr]^{1\times \mu} \ar[d]_{1\times 1\times f} & &
G\times \clM \ar[rrr]^{\mu} \ar[d]_{1\times f} & & &
\clM \ar[d]^{f} \\
G\times G\times \clN \ar[rr]^{1\times \nu} \ar@2[rru]^{1\times \sigma} & &
G\times \clN \ar[rrr]^{\nu} \ar@2[rrru]^{\sigma} & & & \clN
}
\end{array}
\qquad
\begin{array}[t]{c}
\\ \\ \\ \mbox{and}
\end{array}
\qquad
\xymatrix@R=1.5pc@C=4pc{
G\times G\times\clM \ar@/^.8pc/[r] \ar@/_.8pc/[r] \ar[d]
& G\times G\times\clN \ar[d] \\
G\times\clM \ar@{.>}@/^.8pc/[r] \ar@/_.8pc/[r] \ar[d]
& G\times\clN \ar[d] \\
\clM \ar@{.>}@/^.8pc/[r] \ar@/_.8pc/[r] & \clN 
}
$$
\par\vspace{4mm}\par
Observe that if we already know that all pairs $(f,\sigma)$ in a diagram are
1-morphisms of $G$-stacks, then we only have to check 2-commutativity of the
"lower stage" of the prisms. \qed
\end{rema*}


The 2-category $G\tir\frSt/S$ has arbitrary projective and inductive limits.
In particular $G\tir\frSt/S$ has fibred products, defined in the obvious way, so
we have the notion of a \emph{2-cartesian square}.
Any stack $\clM$ over $S$ gives a trivial $G$-stack $(\clM,\pr_2)$
and this gives a 2-functor $\imath:\frSt\to G\tir\frSt$. The invariants and
coinvariants are the 2-adjoints of this functor~:

\begin{defi*} \label{defiFxPtsQuot}
Let $G$ be a sheaf in groups over $S$ and $\clM$ a $G$-stack over $S$.
\begin{trivlist}
\itemn{i} A \emph{stack of fixed points} $\clM^G$ is a stack that
represents the 2-functor $\frSt^\circ\to\frCat$ defined by
$$
F(\clN)=\clH om_{G\tir\frSt}(\imath(\clN),\clM) 
$$
(the latter is the stack of~\ref{defiActionStack}(iii), and $\frCat$ is the
2-category of categories).
\itemn{ii} A \emph{quotient stack} $\clM/G$ is a stack that
represents the 2-functor $\frSt\to\frCat$ defined by
$$
F(\clN)=\clH om_{G\tir\frSt}(\clM,\imath(\clN)) 
$$
\end{trivlist}
\end{defi*}

\begin{rema*} \label{remaTrivSt}
There are in fact several candidates for the notion of a trivial action,
needed to define fixed points and quotients.
The trivial $G$-stacks as defined in \ref{defiActionStack}(iii) form
a full 2-subcategory of $G\tir\frSt$, denoted $\frT\tir G\tir\frSt$.
Its essential image in $G\tir\frSt$ defines the \emph{essentially trivial}
$G$-stacks and we will denote it by $\frE\frT\tir G\tir\frSt$.
The final picture of the factorization of $\imath$  is~:
$$
\xymatrix{
\frSt \ar[r]_>>>>>{
\begin{array}{c}
\uparrow \\
{}_{{\rm ess.} \: {\rm surjective}} \\
{}^{\mathbf{not} \: {\rm f.} \: {\rm faithful}} \\
\end{array}}
& \frT\tir G\tir\frSt \ar[r]^>>>>>{\sim}_>>>>>{
\begin{array}{c}
\uparrow \\
{}_{{\rm eq. } \: {\rm of}} \\
{}^{\rm categories} \\
\end{array}}
& \frE\frT\tir G\tir\frSt \; \ar@{^(->}[r] & G\tir\frSt}
$$
In definition \ref{defiFxPtsQuot} we could have chosen $\clN$
among either of these categories of trivial $G$-stacks. Here the
crucial point is to note the
uncommon feature of $G\tir\frSt$ that, unlike most usual categories
where we can consider group actions (e.g. sets, modules, algebras,
varieties, schemes,...), the quotient of an object with trivial $G$-action
is \emph{not} the object itself. The example of $\clM=S$ and $G$ acting
trivially is in the mind of everyone~: then the stack-quotient is $[S/G]=BG$,
and more generally, for any $\clM$ with trivial $G$-action we should have
$\clM/G=\clM\times BG$.
For this reason it would be meaningless to choose trivial objects among
a full subcategory "$?$" of $G\tir\frSt$, because then we would have
$$
\Hom_{G\tir\frSt}(S,\clN)=\Hom_{\:?}(BG,\clN)=\Hom_{G\tir\frSt}(BG,\clN)
=\Hom_{\:?}(BG\times BG,\clN)=\dots
$$
and so on. This is why $\frSt$ is the only possibility. \qed
\end{rema*}

\begin{prop*} \label{propFxPtSt}
Let $G$ be a sheaf in groups over $S$ and $\clM$ a $G$-stack over $S$.
Then there exists a stack of fixed points $\clM^G$, and its formation
commutes with base change on $S$. The essentially
trivial stacks of remark~\ref{remaTrivSt} are the stacks
$(\clM,\mu)$ isomorphic to $\imath(\clM^G)$.
\end{prop*}

\begin{proo}
From the definition we must have
$\clH om_{\frSt}(\clN,\clM^G)=\clH om_{G\tir\frSt}(\imath(\clN),\clM)$.
From the particular case $\clN=S$ we deduce
$\clM^G=\clH om_{\frSt}(S,\clM^G)=\clH om_{G\tir\frSt}(\imath(S),\clM)$.
This is the stack of $G$-invariant sections of $\clM$, whose objects over a base
$T$ are pairs $(x,\{\alpha_g\}_{g\in G(T)})$ where $x\in \clM(T)$ and
$\alpha_g:x\to g.x$ are isomorphisms such that
$g.\alpha_h\circ \alpha_g=\alpha_{gh}$ for all sections $g,h\in G(T)$.
The second assertion follows because, by definition, an essentially trivial
$G$-stack is a $G$-stack $(\clM,\mu)$ such that there exists an isomorphism
$(\clM,\mu)\simeq \imath(\clN)$ for some $\clN$. Taking fixed points
and then $\imath$, we obtain the result.
\end{proo}

\begin{prop*} \label{propQuotSt}
Let $G$ be a sheaf in groups over $S$ and $\clM$ a $G$-stack over $S$.
Then there exists a quotient stack $\clM/G$, and its formation commutes
with base change on $S$.
\end{prop*}

\begin{proo}
We define a prestack $\clP$ as follows~: sections of $\clP(T)$
are sections of $\clM(T)$, and morphisms in $\clP(T)$ between
$x$ and $y$ are pairs $(g,\varphi)$ with $g\in G(T)$ and $\varphi:g.x\to y$
a morphism in $\clM(T)$. Let $\clM/G$ be the stack associated to $\clP$.
It is straightforward to check the universal 2-property.
\end{proo}

\section{Group actions on algebraic stacks}

In this section we will prove algebraicity of fixed points and
quotients for certain algebraic groups $G$ acting on algebraic stacks.
We consider the category $G\tir\frAlgSt/S$ of algebraic $G$-stacks
over $S$~: this is defined to be the full subcategory of $G\tir\frSt/S$
of $G$-stacks whose underlying stack is algebraic. In particular all
definitions of~\ref{defiActionStack} apply, so we do not have to rewrite
them. The definitions of~\ref{defiFxPtsQuot} carry on in an obvious way,
namely the \emph{algebraic stack of fixed points} represents a 2-functor
$\frAlgSt^\circ\to\frCat$, and the \emph{quotient algebraic stack}
represents a 2-functor $\frAlgSt\to\frCat$.

Before we go further we recall a few examples~:

\begin{exams*} \label{examplesGstacks}
\begin{trivlist}
\itemn{i} Let $G$ be a flat, separated group scheme of finite presentation
over $S$. Then the sheaf $\Aut(G)$ acts on the stack of
$G$-torsors $BG$ by twisting the action~: given $\theta\in\Aut(G)$
and $E\to T$ a $G$-torsor over $T/S$, the twisted action is defined
by $g*e=\theta(g).e$.
\itemn{ii} Let $\overline{\clM}_{g,n}$ be the stack of stable curves
of genus $g$  with $n$ marked points. Then the symmetric group $\frS_n$
acts on it by permuting the marked points.
\itemn{iii} Let $\clM_g(n)$ be the stack of smooth curves of genus $g$
together with a level $n$ structure, i.e. an isomorphism
$\varphi:C[n]\isomto (\zmod{n})^{2g}$. Then $G=\GL_{2g}(\zmod{n})$
acts on $\clM_g(n)$ by twisting the level structure.
\itemn{iv} Let $\clX_1(N)$ be the stack of elliptic curves together
with a "point of $N$-torsion" (see [KaMa]). Then $G=(\zmod{n})^\times$
acts on $\clX_1(N)$ by acting on the point of $N$-torsion.
\end{trivlist}
\end{exams*}

\subsection{Fixed points}

\begin{theo*} \label{theoFPproper}
Let $G$ be a proper, flat group scheme of finite presentation over $S$.
Let $\clM$ be an algebraic $G$-stack, with diagonal locally of finite
presentation over $S$. Then the fixed point stack $\clM^G$
(prop.~\ref{propFxPtSt}) is algebraic (so it is a fixed point stack
in $\frAlgSt$). The morphism $\epsilon:\clM^G\to \clM$ is representable,
separated and locally of finite presentation. The formation of $\clM^G$
commutes with base change on $S$.
\end{theo*}

\begin{proo}
It is enough to show that the morphism $\clM^G\to \clM$ is representable
with the desired properties. So let $f:T\to \clM$ be a 1-morphism, corresponding
to an object $x\in \clM(T)$. The fibre product $\clM^G\times_{\clM} T$ is
the sheaf whose sections over $T'/T$ are collections of isomorphisms
$\{\alpha_g:x\simeq g.x\}_{g\in G(T')}$ such that for all sections $g,h\in G(T')$
we have $g.\alpha_h\circ \alpha_g=\alpha_{gh}$. Denote by $x_1$ and
$x_2$ the objects of $\clM(G\times T)$ corresponding to the 1-morphisms
$\pr_2\circ (\Id_G\times f)$ and $\mu\circ (\Id_G\times f)$. Reformulating
what we said above, there is a closed immersion, locally of finite presentation,
$\clM^G\times_{\clM} T \hookrightarrow \Hom_T(G_T,\Isom_{G_T}(x_1,x_2))$.
With our assumptions, a section of this $\Hom$ sheaf gives, via its graph, a closed
subspace of $G_T\times \Isom_{G_T}(x_1,x_2)$ which is proper and flat over
the base, being isomorphic to $G_T$ via the first projection. So the sheaf
is representable by the corresponding open constructible subspace of the
Hilbert space, which is algebraic, separated and locally of finite presentation
by (Artin \cite{Ar} cor. 6.2).
The result follows.
\end{proo}

\begin{remas*}
\begin{trivlist}
\itemn{i}
If $\clM$ is representable, then $\clM^G$ is representable also,
and so the fixed points of $\clM$ as a space or as a stack are the same
(in general the Yoneda functor from spaces into stacks commutes with
projective limits when they exist,
but not with inductive limits; see also~\ref{remaQuotDependsSpSt}).
\itemn{ii}
For an essentially trivial $G$-stack $(\clN,\nu)$ (see \ref{remaTrivSt} and
\ref{propFxPtSt}), arbitrary $G$-morphisms $(f,\sigma):\clN\to\clM$ still
factor through $\clM^G$, because $(\clN,\nu)\simeq (\clN^G,\pr_2)$. This
factorization is of course not unique.
\end{trivlist}
\end{remas*}

If we relax the assumption of properness of $G$, it does not seem plausible
that we can say much on representability of $\epsilon:\clM^G\to \clM$,
at least if the diagonal of $\clM$ is not flat. If we put conditions on $G$
and on the diagonal of $\clM$ such as flatness or smoothness, then it may be
that using arguments such as these developped in \cite{SGA3}, tome 2, we obtain
representability of $\clM^G\to \clM$ in some cases.

We now wish to give more properties of the morphism $\epsilon:\clM^G\to \clM$
when the group $G$ is \emph{finite} and $\clM$ is separated. In this case, it
is not possible to deduce that the space $\Hom_T(G_T,\Isom_{G_T}(x_1,x_2))$ in
the proof above is finite, or even proper, because $G$ may be ramified. However
we will deduce the corresponding property for $\epsilon$ by giving a slightly
different construction of $\clM^G$. We start with a lemma~:

\begin{lemm*}
Let $Q$ be a finite flat scheme of finite presentation over $S$.
Let $\clM$ be an algebraic stack locally of finite presentation over $S$.
Then the stack $\clH om_S(Q,\clM)$ of morphisms of stacks from $Q$ to
$\clM$ is algebraic and locally of finite presentation over $S$.
\end{lemm*}

\begin{proo}
Let's note $\clH:=\clH om_S(Q,\clM)$ and $n=[Q:S]$. Notice that, given
an $S$-scheme $T$, we have $\clH(T)=\clM(Q\times T)$. From this and the fact
that $Q$ is affine, after algebraicity is proved it will follow that $\clH$
is locally of finite presentation over $S$ because given a filtering inductive
system of $S$-algebras $A_i$, we have isomorphisms
$$
\underset{\longto}{\lim}\,\clH(A_i)\simeq
\underset{\longto}{\lim}\,\clM(\clO_Q\otimes A_i)\simeq
\clM(\underset{\longto}{\lim}\,\clO_Q\otimes A_i)\simeq
\clM(\clO_Q\otimes \underset{\longto}{\lim}\,A_i)\simeq
\clH(\underset{\longto}{\lim}\,A_i)
$$
Now we show that the diagonal of $\clH$ is representable, separated and quasi-compact.
It is enough to study the sheaf $\Isom_{\clH_T}(x,y)$ for two fixed objects $x,y\in \clH(T)$.
These correspond to objects $\eta\in\clM(Q\times T)$ and $\xi\in\clM(Q\times T)$, and
$$
\Isom_{\clH_T}(x,y)=\Hom_T\big(Q_T,\Isom_{\clM_{Q\times T}}(\eta,\xi)\big)
$$
Here the sheaf $I:=\Isom_{\clM_{Q\times T}}(\eta,\xi)$ is representable and of
finite presentation over $Q_T$ (it is locally of finite presentation because
$\clM$ is, by \cite{EGA}, I, 6.2.6. which extends to stacks). It keeps these
properties as a $T$-sheaf. Let us introduce the functor $H_n$ which is the
component of the full Hilbert functor of $Q_T\times I$ parametrizing 0-dimensional
subspaces of length $n$. It is representable by a separated algebraic space locally
of finite presentation (Artin \cite{Ar} cor. 6.2), and in fact the length $n$
component is quasi-compact because $Q_T\times I$ is. Now, the graph of a morphism
$Q_T\to I$ defines a point in $H_n$ (by separation of $I$), such that the restriction
of the first projection $Q_T\times I\to Q_T$ is an isomorphism. The sheaf
$\Isom_{\clH_T}(x,y)$ is thus isomorphic to the corresponding constructible
open subspace of $H_n$.
By constructibility this open immersion is quasi-compact (\cite{EGA},
$0_{\rm III}$, 9.1.5), and, of course, separated.

Now let $U\to \clM$ be an atlas~; we can choose $U$ separated.
Then I claim that $V:=\Hom_S(Q,U)$ will be an atlas for $\clH$.
First, by Artin's result again $V$ is representable and locally
of finite presentation. As $\clH$ is also locally of finite presentation
this shows that the map $V\to \clH$ has the same property.
Thus we only have to prove that it is formally smooth and
surjective.

To prove surjectivity take an algebraically closed field $k$ and a morphism
$\Spec(k)\to \clH$ i.e. a morphism $f:Q_k\to \clM_k$. Then $Q_k$ is an artinian
scheme, hence a sum of local artinian $k$-schemes, so we reduce to the local
case. By surjectivity of $U\to \clM$, the image of the underlying point
of $Q_k$ lifts to $U_k$, and by smoothness the whole morphism lifts.

It remains to prove formal smoothness. Let $A\to A_0$ be a surjection of
artinian rings with nilpotent kernel. Assume we have a 2-commutative
diagram
$$
\xymatrix{
V \ar[r] \ar@2+<7mm,-7mm>^>>>\delta & \clH \\
\Spec(A_0) \ar[u] \ar[ru] }
\qquad
\begin{array}[t]{c}
\\ \mbox{meaning that we have}
\end{array}
\qquad
\xymatrix{
U_{A_0} \ar[r] \ar@2+<7mm,-7mm>^>>>\delta & \clM_{A_0} \\
Q_{A_0} \ar[u] \ar[ru] }
$$
As $Q_{A_0}$ is artinian, by smoothness of $U_A\to \clM_A$, the
map $Q_{A_0}\to U_{A_0}\to U_A$ immediately lifts to $Q_A\to U_A$, and we are done.
\end{proo}

\begin{remas*}
\begin{trivlist}
\itemn{i}
If $Q=S[\varepsilon]/\varepsilon^2$ we recover the tangent stack $T(\clM/S)$,
and the lemma gives a proof of its algebraicity which is simpler than in
\cite{LMB}, chap. 17. If $Q=S[x]/x^n$ we get the stack of $n$-truncated
arcs in $\clM$. If $Q$ is sum of $n=[Q:S]$ copies of $S$ then
$\clH om_S(Q,\clM)=\clM^n$ so the result is trivial.
\itemn{ii}
The result of \cite{OS} of representability of Quot functors for Deligne-Mumford
stacks does not allow to derive algebraicity of $\clH om_S(Q,\clM)$
because one can not express this stack as an open substack of the "Hilbert space"
(graph morphisms of algebraic stacks are no longer closed immersions, not even
for $\clM$ separated).
\end{trivlist}
\end{remas*}

\begin{prop*}
Let $G$ be a finite, flat group scheme of finite presentation over $S$.
Let $\clM$ be an algebraic $G$-stack, locally of finite presentation over $S$.
Then the morphism $\clM^G\to \clM$ is furthermore quasicompact, and enjoys
any property enjoyed by the diagonal of $\clM$, by closed immersions,
and stable by composition. In particular it is proper if $\clM$ is separated.
\end{prop*}

\begin{proo}
Throughout, we will omit the description of the morphisms of the different stacks
introduced, since they are obvious and quite lengthy to write completely.
By the lemma applied to $Q=G$ the stack $\clH=\clH om(G,\clM)$ is algebraic.
We now define two morphisms $a,b:\clM \to \clH$. Let $x\in \clM(T)$, corresponding
to a morphism $f:T\to \clM$, and look at the compositions
$$
\xymatrix{
G\times T \ar[r]^{\Id\times f} & G\times\clM \ar@<-.4ex>[r]_>>>>>{\pr_2}
\ar@<.4ex>[r]^>>>>>{\mu} & \clM}
$$
Then we define $a(x)=(\mu\circ (\Id_G\times f))^*(x)$
and $b(x)=(\pr_2\circ (\Id_G\times f))^*(x)=x_{G_T}$. In more naive terms,
$a(x)=(g\mapsto g.x)$ and $b(x)=(g\mapsto x)$. Now look at the fibre product
defined by the diagram
$$
\xymatrix{
\clN \ar[r] \ar[d] & \clM \ar[d]^b \\
\clM \ar[r]^>>>>>a & \clH
}
$$
An object of
$\clN$ is a pair $(x,\psi^x:a(x)\simeq b(x))$ where $\psi^x$ consists
in isomorphisms $\psi^x_g:g.x\simeq x$. We define a closed substack
$\clZ\subset \clH$ by considering the morphisms $\psi:G\to \clM$
such that for all sections $g,h\in G(T)$ we have $g(\psi_h)\circ \psi_g=\psi_{gh}$.
The stack $\clM^G$ is isomorphic to the preimage of $\clZ$ in $\clN$.
The morphism $\epsilon:\clM^G\to \clM$ is the first projection.
Finally, it is not hard to check that $\clM^G$ is locally of finite
presentation, using that it is the case for $a$ and for $\clH$
and its diagonal.


It remains to prove the properties of the morphism $\clM^G\to \clM$.
First we look at the morphism $b:\clM\to \clH$. Let
$U\to \clH$ be a morphism, corresponding to an object
$\xi\in \clM(G\times U)$. The fibre product $\clM\times_{\clH} U$
is the stack of triples $(T,\eta,\alpha)$ composed of a map of schemes
$T\to U$, an object $\eta\in \clM(T)$ and an isomorphism $\alpha$ between
$\eta_{G_T}$ and $\xi_{G_T}$. By fppf descent, this is none other than
the functor of descent data for $\xi$ with respect to the fppf covering
$G_U\to U$. It is represented by a closed sub-algebraic space of
$\Isom_{G_U\times_U G_U}(\pr_1^*\xi,\pr_2^*\xi)$. This space inherits
the properties such as quasi-compactness and separatedness
of the diagonal of $\clM$. It follows that $b$ has these properties,
and similarly for $\clN$ and $\clM^G$.
\end{proo}

In the case of a finite constant group $G$, everything is much simpler
because we have $\clH om_S(G,\clM)=\clM^n$ as noticed earlier, and we do
not need to make assumptions on $\clM$.

\begin{exam*}
The following example shows that the morphism $\clM^G\to \clM$ needs not
be a monomorphism of algebraic stacks, although $\imath(\clM^G)\to \clM$
is necessarily a monomorphism of $G$-algebraic stacks, because of the
2-universal property.
Let $\clM_{g,2}$ be the stack of smooth 2-pointed curves of genus $g$
(see~\ref{examplesGstacks}(ii)).
It has an action of the symmetric group $\frS_2$. Let $(C,a,b)$ be a curve
over a base $S$, and suppose that $C$ has two distinct automorphisms
$\sigma_1$ and $\sigma_2$ which exchange the marked points. Then these
give two morphisms $S\to (\clM_{g,2})^{\frS_2}$, and the compositions
$S\to \clM_{g,2}$ are equal as morphisms of algebraic stacks. However,
they are not equal as morphisms of $\frS_2$-algebraic stacks because the
maps $\sigma_1,\sigma_2$ enter in the definition of such a morphism.
\end{exam*}

\begin{exam*}
The following example shows that "fixed points" and
"coarse moduli space" do not commute.
Let $Q=\{\pm 1,\pm i,\pm j,\pm k\}$ be the quaternion group, of order 8.
Its unique involution generates its center $Z$, and $G=Q/Z\simeq
\zmod{2}\times \zmod{2}$ is not isomorphic to a subgroup of $Q$.
There is a faithful action of $G$ on $Q$ by conjugation, whence an
action of $G$ on $BQ$ (see~\ref{examplesGstacks}(i)). For the trivial
$Q$-torsor $x:Q\to S$, for all $g$ the left multiplication by $g^{-1}$
is an isomorphism $g.x\simeq x$, and however there is no $G$-linearization
so $x$ is not a fixed point. Actually $(BQ)^G$ is empty, whereas
the moduli space of $BQ$ is $S$ and we have $S^G=S$ for the induced
action.
\end{exam*}

\begin{exam*}
The following example shows that $\clM^G$ may not be algebraic when $G$ is
not proper. If $H$ is a commutative group scheme and $G$ a group scheme
acting trivially on $BH$, then an objet of $(BH)^G$ is an $H$-torsor $x$
together with a morphism $G\to \Aut(x)=H$, so $(BH)^G=BH\times\Hom(G,H)$.
This stack is not algebraic in general, though for special groups $G,H$
it may be the case (for instance if both $G,H$ are of multiplicative type
--- see \cite{SGA3}, tome 2 again).
\end{exam*}


\subsection{Quotients}

Let $G$ be a flat, separated group scheme of finite presentation over $S$.
By a $G$-torsor over an $S$-scheme~$T$, we will mean an algebraic space
with $G$-action $p:E\to T$ that locally on $T$ is isomorphic to the trivial
$G$-space $G\times T$. In general such a torsor will not be a scheme, unless
if for example $G$ is quasi-affine.

Let $\clM$ be a $G$-algebraic stack over~$S$. In case
$\clM=X$ is an algebraic space, the quotient of~\ref{propQuotSt} is known
under the more familiar decription of the stack of $G$-torsors with an
equivariant morphism to $X$. It is traditionnally denoted $[X/G]$, to avoid
confusion with a hypothetical \emph{quotient algebraic space}, but when
$\clM$ is a general stack no such confusion is possible so it is natural
to suppress the brackets.

For general $\clM$  we can still define a stack whose objects are
$G$-torsors $p:E\to T$ with an equivariant morphism $(f,\sigma):E\to \clM$.
More precisely we define a stack $(\clM/G)^*$ whose sections over $T$ are
triples $t=(p,f,\sigma)$ as above, and the isomorphisms between $t$ and $t'$
in $(\clM/G)^*$ are pairs $(u,\alpha)$ with a $G$-morphism $u:E\to E'$ and
a 2-commutative diagram of $G$-stacks (see~\ref{remaCommDiag})
$$
\xymatrix@C=1pc{
E \ar[rd]_{(f,\sigma)} \ar[rr]^u & & E' \ar[ld]^{(f',\sigma')}
\ar@2{<=}+<-12mm,-7mm>_(.6)\alpha \\
& \clM &
}
$$

\begin{theo*}
Let $G$ be a flat, separated group scheme of finite presentation over $S$.
Let $\clM$ be a $G$-algebraic stack over $S$. Then the quotient stack $\clM/G$
(prop.~\ref{propQuotSt}) is isomorphic to the stack of $G$-torsors
$(\clM/G)^*$, and it is algebraic (so it is a quotient stack in $\frAlgSt$).
The canonical morphism $\pi:\clM\to \clM/G$ is the universal torsor over
$\clM/G$. The formation of $\clM/G$ commutes with base change on $S$.
\end{theo*}

\begin{proo}
There are two things to show.
First, we explain why $\clM/G\simeq (\clM/G)^*$. Let $\clM/G$
be the quotient as described in \ref{propQuotSt}, which is the stack associated
to a prestack $\clP$. We define a morphism $u:\clP\to (\clM/G)^*$ by sending
an object $x\in \clM(T)$ viewed as a map $x:T\to \clM$, to the trivial torsor
together with the map $G\times T\to \clM$ given by $\mu\circ (\Id\times x)$,
which is clearly equivariant. The image of a morphism $(g,\varphi):x\to y$ in
$\clP$ is the multiplication by $g$ (as a map of torsors).
This morphism $u$ extends to a morphism of stacks $u':\clM/G\to (\clM/G)^*$.
It is clearly fully faithful, and also locally essentially surjective by the
definition of a torsor. So it is an isomorphism of stacks. From now on we
identify $\clM/G$ and $(\clM/G)^*$.

Second, we prove algebraicity. We keep the above notations of $t=(p,f,\sigma)$
for sections of $\clM/G$ and $\varphi=(u,\alpha)$ for morphisms between
$t$ and $t'$. Note that there is a morphism $\omega:\clM\to BG$ obtained
by forgetting the maps to $\clM$. To study the diagonal of $\clM/G$, we
take $t,t'\in (\clM/G)(T)$, then $\omega$ induces a morphism
$\Isom_T(t,t')\to \Isom_{BG}(E,E')$ given by $(u,\alpha)\mapsto u$.
The latter space is algebraic, and the fibre of this projection above
an isomorphism $u:E\to E'$ is the closed (algebraic) subspace of
$\Isom_{\clM_T}(E,u^*E')$ of 2-$G$-isomorphisms. This shows that
$\Isom_T(t,t')$ is representable, separated, quasi-compact. From the
fact the morphism $S\to BG$ is fppf, and the obvious 2-cartesian diagram
$$
\xymatrix{
\clM \ar[r] \ar[d] & S \ar[d] \\
\clM/G \ar[r] & BG
}
$$
we deduce that $\clM\to \clM/G$ is fppf, and by composition with an atlas
of $\clM$ we get an fppf presentation of $\clM/G$, whence the result.
\end{proo}

\begin{rema*} \label{remaQuotDependsSpSt}
It is clear that if $\clM$ is representable, then the quotient $\clM/G$
depends on if we compute it in the category of spaces or of stacks. Any
algebraic space $X$ with non-free action of a finite group $G$ has a
quotient space $X/G$, distinct from the quotient stack.
\end{rema*}

\begin{exam*}
Let $\clM$ be a $G$-algebraic stack over $S$, so we have morphisms
$
\xymatrix{G\times\clM \ar@<-.4ex>[r]
\ar@<.4ex>[r]^>>>>>{\mu,\pr_2} & \clM}
$.
Given a sheaf $F$ on the smooth-\'etale site of $\clM$, a $G$-linearization
of $F$ is an isomorphism $\alpha:\mu^*F\simeq \pr_2^*F$ which is compatible
with associativity~: $(m\times \Id_{\clM})^*\alpha=(\Id_G\times \mu)^*\alpha$.
We define a (smooth-\'etale) \emph{$G$-sheaf} on $\clM$
to be a pair $(F,\alpha)$ as above.
We can look at the stack of
invertible $G$-sheaves (with obvious isomorphisms of $G$-sheaves between them),
denoted $\clP ic^G(\clM)$, and it is easy to see that we have canonical
isomorphisms of stacks $\clP ic (\clM/G)\simeq \clP ic^G(\clM)\simeq \clP ic(\clM)^G$.
In particular if $\clP ic(\clM)$ is algebraic and $G$ is proper, flat,
of finite presentation, we obtain algebraicity of the first two stacks,
by theorem~\ref{theoFPproper}.
\end{exam*}

\end{document}